\documentstyle[12pt]{article}
\input amssym.def
\newtheorem{th}{Theorem}[section]
\newtheorem{df}[th]{Definition}
\newtheorem{prop}[th]{Proposition}

\newsymbol\rtimes 226F
\newsymbol\ltimes 226E
\newcommand{\mod}{{{\rm mod\,}}}
\newcommand{\tN}{{{\tilde N}}}
\newcommand{\tM}{{{\tilde M}}}
\newcommand{\talpha}{{\tilde \alpha}}
\newcommand{\tbeta}{{\tilde \beta}}
\newcommand{\tr}{{\rm tr}}

\begin{document}
\title{Classification of actions of
       discrete \\amenable groups on 
      \\strongly amenable subfactors of type III$_\lambda$}
\author{Toshihiko Masuda\\
Department of Mathematical Sciences,\\
University of Tokyo, Komaba, Tokyo, 153. JAPAN}
\date{}
\maketitle
\begin{abstract}
Using the continuous decomposition, we classify strongly free actions of
discrete amenable groups on strongly amenable subfactors of type III$_\lambda$,
$0<\lambda<1$. Winsl\o w's fundamental homomorphism is a 
complete invariant. This removes the extra assumptions in the classification
theorems of Loi and Winsl\o w and gives a complete classification up to 
cocycle conjugacy.
\end{abstract}


\section{Introduction}
In the theory of operator algebras, the study of automorphisms is one of 
the most important topics. Especially since Connes's work \cite{C2},
much progress has been made on the classification of the 
actions of discrete amenable groups on injective factors.

In the subfactor theory, various studies of automorphisms have been done. 
In \cite{P2}, Popa has introduced the notion of proper outerness of 
automorphisms
and proved that the strongly outer actions of discrete amenable groups 
on strongly amenable subfactors of type II$_1$ are classified by
the Loi invariant (See \cite{L1}). 
(In \cite{CK}, Choda and Kosaki have introduced the same property 
independently and they call it strong outerness.)

In the case of subfactors of type III$_\lambda$, $(0<\lambda<1)$, 
partial results on classification of group actions have been obtained 
by Winsl\o w and Loi. (\cite{L2}, \cite{W1}, \cite{W2}.)
In \cite{W1} and \cite{W2}, Winsl\o w has introduced the strong freeness
and the fundamental homomorphism for actions. 
He has classified the strongly free actions of discrete amenable groups on 
subfactors of type III$_\lambda$ for groups having the character
lifting property. His fundamental homomorphism is a complete invariant.
In \cite{L2}, Loi gave a classification theorem when $G$ is finite. Their 
idea of the proof is 
that they reduce the classification problem to the type II$_\infty$ case 
using the discrete decomposition and apply Popa's classification result. 
The most difficult points of their proofs are to reduce the problem to
the type II$_\infty$ case. Because of this difficulty, they made extra 
assumptions such as the character lifting property or finiteness for groups.  
But it seems difficult to generalize their method to the arbitrary 
discrete amenable group case.    

Our idea of a proof is using the continuous decomposition instead of the 
discrete decomposition based on the method in \cite{ST1}, \cite{ST2}. 
But in this case, we treat only factors of type III$_\lambda$, $0<\lambda<1$, 
with the trivial characteristic invariants, so the proof is less complicated 
than those in \cite{ST1} and \cite{ST2}. By using the continuous decomposition,
we can more easily reduce the classification probelm to the type II$_\infty$ 
case than using the discrete decomposition and this method is valid for
arbitrary discrete amenable groups. \\
\vspace{3pt}\\
{\bf Acknowledgement.} The author is grateful to Prof. M. Izumi for proposing 
this problem to him and helpful suggestions, Prof. Y. Kawahigashi for 
fruitful comments and constant encouragement.   

\section{Preliminaries}
  In this section, we recall several results about group actions on subfactors,
and fix notations. The facts stated in this section are found in 
\cite{CK}, \cite{L1}, \cite{L2}, \cite{P2}, \cite{W1}, \cite{W2}, \cite{W3}.
 
  Let $N\subset M$ be an inclusion of factors with finite index and 
$N\subset M \subset M_1 \subset M_2 \subset \cdots $ the Jones tower.   
(Throughout this paper, we always assume that conditional expectations 
are minimal in the sense of \cite{H} and inclusions of factors of type II are 
extremal. ) For 
$\alpha \in {\rm Aut}(M,N)$, we extend $\alpha$ to $M_k$ such that 
$\alpha(e_k)=e_k$ inductively, where $e_k$ denotes the Jones projection for 
$M_{k-1}\subset M_k$. 

First we recall the Loi invariant and the strong outerness of 
group actions.

\begin{df}[{\cite[Section 5]{L1}}]
With above notations, Put
$$\Phi(\alpha):=\{\alpha|_{M'\cap M_k}\}_k.$$
We call $\Phi$ the Loi invariant for $\alpha$.
\end{df}

\begin{df}[{\cite[Definition 1]{CK}}, {\cite[Definition 1.5.1]{P2}}]
An automorphism \\
$\alpha\in{\rm Aut}(M,N)$ is said to be properly outer or
strongly outer if we have no non-zero $a\in\bigcup_k M_k$ satisfying 
$\alpha(x)a=ax$ for all $x\in M$. The action $\alpha$ of $G$
 on $N\subset M$ is said to be strongly outer if 
$\alpha_g$ is strongly outer except for $g=e$.    
\end{df}

The most important result on classification of actions of groups 
on subfactors has been obtained by Popa.
\begin{th}[{\cite[Theorem 3.1]{P2}}]
Let $N\subset M$ be a strongly amenable inclusion of factor of type $II_1$ 
and $G$ a countable discrete amenable group.

If $\alpha$ and $\beta$ are strongly outer actions of $G$ on $N\subset M$,
then $\alpha$ and $\beta$ are cocycle conjugate if and only if 
$\Phi(\alpha)=\Phi(\beta)$.  
\end{th}

For type II$_\infty$ inclusions, we have the following result due to
Popa and Winsl\o w. 

\begin{th}[{\cite[Theorem 2.1]{P2}}, {\cite[Theorem 4.3]{W1}}]\label{th:infty} 
Let $N\subset M$ a strongly ame-
nable inclusion of factors of type II$_\infty$.
If $\alpha$ and $\beta$ are actions of countable discrete amenable group $G$
on $N \subset M$, then $\alpha$ and $\beta$ are cocycle conjugate if and only 
if $\Phi(\alpha)=\Phi(\beta)$ and $\mod(\alpha)=\mod(\beta)$.  
\end{th}

Let $N\subset M$ be an arbitrary inclusion of factors with the common flow of 
weights. Fix a normal state of $N$ and take a crossed product of $N\subset M$
by the modular automorphism. 
Put ${\tilde N} \subset {\tilde M}:=N\rtimes_{\sigma^\phi}{\bf R} \subset
M\rtimes_{\sigma^{\phi\circ E}}{\bf R}$, where $E$ is the minimal conditional
expectation from $M$ onto $N$.
Let   
${\tilde \alpha}$ be the canonical extension of $\alpha$ to 
${\tilde N}\subset{\tilde M}$ (\cite{CT}, \cite{HS}), i.e., 
$${\tilde \alpha}(x):=\alpha(x), \quad x\in M,$$
$${\tilde \alpha}(\lambda(t)):=(D\phi\,\alpha^{-1}:D\phi)_t\,\lambda(t),$$
where $\lambda(t)$ is the usual implementing unitary.
The notions of strong freeness for automorphisms and the fundamental
homomorphism are introduced by Winsl\o w in \cite{W1}, \cite{W2}.

\begin{df}[{\cite[Definition 3.2]{W1}}, {\cite[Definition 4.2]{W2}}]
An automorphism \\
$\alpha\in{\rm Aut}(M,N)$ is said to be strongly free
if we have no non-zero $a\in\bigcup_k {\tilde M_k}$ satisfying 
${\tilde \alpha}(x)a=ax$ for all $x\in {\tilde M}$.
For an action $\alpha$ of G on $N\subset M$ is said to be strongly free if 
$\alpha_g$ is strongly free except for $g=e$.    
\end{df}

According to \cite{W1} and \cite{W2}, we set 
$$\Upsilon(\alpha):=\{{\tilde \alpha}|_{{\tilde M}'\cap{\tilde M_k}}\}_k$$
and we call this the fundamental homomorphism.


\section{Classification of actions}
Throughout this section, we assume that inclusions of factors of type 
III$_\lambda$
are strongly amenable in the sence of Popa. (See \cite{P1} and \cite{P2}.)

The following theorem is the main result of this paper. 
\begin{th}\label{main}
Let $N\subset M$ be a strongly amenable inclusion of factors of 
type III$_\lambda$, $0<\lambda<1$, with the common flow of weights. 
Let $G$ be a countable discrete amenable group,
and $\alpha$ and $\beta$ strongly free actions of $G$ on 
$N\subset M$. Then $\alpha$ and $\beta$ are cocycle conjugate
if and only if ${\rm \Upsilon}(\alpha)={\rm \Upsilon}(\beta)$.   
\end{th}
Our idea of proof is that we lift actions  
to inclusions of type II$_\infty$ von Neumann algebras 
using continuous decomposition and apply Popa's result.

The ``only if'' part is obvious, so we only prove the ``if'' part. 
Let $(X, F_t)$ be the flow of weights of $M$. Since $M$ is of 
type III$_\lambda$, $(X, F_t)$ is of the form 
$([0,-\log\lambda), {\rm translation})$.
And we have an isomorphism 
$$(\tN \subset \tM \subset \tM_1 \subset \cdots )\cong
(   L^\infty(X)\otimes Q \subset L^\infty(X)\otimes P \subset
 L^\infty(X)\otimes P_1 \subset \cdots) , $$
where $Q\subset P\subset P_1 \subset \cdots $ is a tower of factors of
type II$_\infty$ and $Q\subset P$ is strongly amenable by assumption.
 
Let $\theta_t$ be the usual trace scaling action of {\bf R} on 
$\tN\subset\tM$. Since ${\tilde \alpha}_g, g\in G$, commutes with 
$\theta_t$, we can consider 
an  action of $G\times {\bf R}$ by setting
$(g,t)\to {\tilde \alpha}_g\theta_t$. If no confusion arises, we also denote
this action of $G\times {\bf R}$ by $\talpha$.
If we prove that two actions of $G\times {\bf R}$, $\talpha$ and $\tbeta$
are cocycle conjugate,  
the proof of \cite[Proposition 1.1]{ST2} also works in this case and we can 
deduce that the canonical extensions of $\talpha$ and $\tbeta$ on 
$N\rtimes_{\sigma^\phi}{\bf R}\rtimes_\theta{\bf R} \subset 
M\rtimes_{\sigma^{\phi\circ E}}{\bf R}\rtimes_\theta{\bf R}
\cong N\otimes B(L^2({\bf R}))\subset M\otimes B(L^2({\bf R}))$ are also 
cocycle conjugate and get the conclusion that $\alpha$ and $\beta$
are cocycle conjugate, since $N$ and $M$ are properly infinite.
 
So our purpose is the classification of actions of 
$G\times {\bf R}$ on $\tN \subset \tM$. 
Note that for $g\in G$, the equality $\tr_\tM\,\talpha_g=\tr_\tM$ holds.

Put $H:=G\times {\bf R}$ and we consider the action $\talpha$ of the 
groupoid $H\ltimes X$ on $Q\subset P$ by the equality 
$$\talpha_g(a):=\int^\oplus_X\,\talpha_{(g,g^{-1}x)}(a(g^{-1}x))\,dx.
   \quad\mbox{(See \cite[Proposition 1.2]{ST2}.)}$$

Put $x_0:=0\in X$ and $H_0:=\{(g,x_0)\in H\ltimes X\, |\, gx_0=x_0\}$.
Then $H_0$ is a discrete amenable group acting on $Q\subset P$.
For $x\in X=[0,-\log\lambda)$, we define $h(x): X\to X$ by
$h(x)y:= y+x$, where sum is taken modulo $-\log\lambda$.
Especially $h(x)x_0=x$.

Here we have the following proposition.

\begin{prop}\label{prop:outer} 
If the action of $G$ is strongly free, then the action of $H_0$ is
strongly outer. 
\end{prop}
{\bf Proof.} 
Assume that the action of $H_0$ is not strongly outer. Then there exists 
$g\in H\backslash \{e\}$ and a nonzero $a\in P_k$ for some $k$ such that 
for every $b\in Q$, we have $\talpha_g(b)a=ab$. Since $\talpha_g$ is not strongly outer and $\mod\talpha_g=1$, we know that $g$ is in $G$.

Set 
$${\tilde a}:=\int^\oplus_X \talpha_{(h(x), x_0)}(a) dx.$$ 
Then an easy computation shows that the equality 
$\talpha_g(b){\tilde a}={\tilde a}b$ holds for every $b\in Q$ and 
this means that action $\alpha$ is not strongly free. \hfill$\Box$ 
\vspace{7pt}\\
{\bf Proof of Theorem \ref{main}}
Let $\alpha$ and $\beta$ be strongly free actions of $G$ on $N\subset M$ 
such that $\Upsilon(\alpha)=\Upsilon(\beta)$. Then we get two actions 
$\talpha$ and $\tbeta$ of the same groupoid $H\ltimes X$. 
So we get two actions $\talpha$ and $\tbeta$ of a discrete amenable group
$H_0$. 

By Proposition \ref{prop:outer}, both actions are strongly outer and by 
assumption both actions have the same Loi invariant and the same module. 
So there exists an automorphism $\theta \in {\rm Aut}(P,Q)$ and  
$u_g\in Z_{\tbeta}(H_0, U(Q))$ such that 
$$\theta\,\talpha_g\,\theta^{-1}={\rm Ad}u_g\,\tbeta_g, \quad g\in H_0.$$  

Set 
$$ \theta_x:=\tbeta_{(h(x),x_0)}\,\theta\,\talpha^{-1}_{(h(x),x_0)},
\quad x\in X\quad \mbox{and}$$

$$ u_{(g,x)}:=\tbeta_{(h(gx)^{-1}gh(x),x_0)}(u_{h(gx)^{-1}gh(x)}),
\quad (g,x)\in H\ltimes X.$$

Then an easy computation shows that $u_{(g,x)}\in Z_\tbeta(H\ltimes X, U(Q))$
and the equality 
$$\theta_{gx}\,\talpha_{g,x}\,\theta^{-1}_x
={\rm Ad}\,u_{(g,x)}\,\tbeta_{(g,x)}$$
holds.

Put ${\tilde \theta}:=\int_X^{\oplus}\theta_x\,dx$
and ${\tilde u}_g:=\int_X^{\oplus}u_{(g,x)}dx.$
Then we get 
$${\tilde \theta}\,\talpha_g\,{\tilde \theta}^{-1}
={\rm Ad}\,{\tilde u}_g\,\tbeta_g$$
and we get the conclusion. \hfill$\Box$

 
\end{document}